\documentclass[11pt,leqno,twoside]{article}
\usepackage{amsthm}
\usepackage{amsmath}
\usepackage{amssymb}
\usepackage[all]{xy}

\title{Substitution Tilings and Separated Nets with Similarities to the Integer Lattice}
\author{Yaar Solomon}

\newenvironment{diagram*}[1]{\begin{align*}{\xymatrix{#1}} \end{align*}} {}

\newenvironment{diagramtow*}[2]{\begin{align*}{\xymatrix{#1}} \qquad {\xymatrix{#2}} \end{align*}} {}

\newcommand{\ceil}[1]{\left\lceil{#1}\right\rceil}

\newcommand{\minus}{\smallsetminus}
\newcommand{\norm}[1]{\left\|{#1}\right\|}

\newcommand{\absolute}[1] {\left|{#1}\right|}
\newcommand{\N}{{\mathbb{N}}}
\newcommand{\Z}{{\mathbb{Z}}}

\newcommand{\R}{{\mathbb{R}}}

\newcommand {\ignore}[1]  {}

\newtheorem{thm}{Theorem}[section]
\theoremstyle{definition}
\newtheorem{definition}[thm]{Definition}
\newtheorem*{Definition}{Definition}
\theoremstyle{plain}
\newtheorem{lem}[thm]{Lemma}
\newtheorem{prop}[thm]{Proposition}
\newtheorem{cor}[thm]{Corollary}

\newtheorem*{example}{Example}

\theoremstyle{remark}

\begin{document}

\date{18.01.2009}
\maketitle

\begin{abstract}
We show that any primitive substitution tiling of $\R^2$ creates a separated net which is biLipschitz to $\Z^2$. Then we show that if $H$ is a primitive Pisot substitution in $\R^d$, for every separated net $Y$, that corresponds to some tiling $\tau\in X_H$, there exists a bijection $\Phi$ between $Y$ and the integer lattice such that $\sup_{y\in Y}\norm{\Phi(y)-y}<\infty$. As a corollary we get that we have such a $\Phi$ for any separated net that corresponds to a Penrose Tiling. 
The proofs rely on results of Laczkovich, and Burago and Kleiner.

\end{abstract}

\section{Introduction} 

\begin{Definition}\label{s.net+biL}
A set $Y\subseteq\R^d$ is called a \emph{separated net}, or a \emph{Delone set}, if there exist constants $R,r>0$ such that every ball of radius $R$ intersects $Y$ and every ball of radius $r$ contains at most one point of $Y$. 
\end{Definition}

\begin{definition}
Let $Y_1$ and $Y_2$ be separated nets. We say that a mapping $\Phi:Y_1\to Y_2$ is \emph{biLipschitz} if there exists a constant $C\ge 1$ such that for every $y,y'\in Y_1$ we have
\[\frac{1}{C}\cdot\norm{y-y'}\le\norm{\Phi(y)-\Phi(y')}\le C\cdot\norm{y-y'},\]  
where $\norm{\cdot}$ is some (any) norm on $\R^d$. $\Phi$ is called a \emph{bounded displacement} if   
\[\sup_{y\in Y_1}\norm{\Phi(y)-y}<\infty.\]
\end{definition}

Consider two equivalence relations on the set of all separated nets. In one, two nets are equivalent if there exists a biLipschitz bijection between them. In the other, the relation holds if there is a bijection which is also a bounded displacement. Since we are dealing with functions between separated nets, it is easy to verify that the second relation refines the first. A natural question is: Is every separated net in $\R^d$ biLipschitz to $\Z^d$? This question was first posed by Gromov in \cite{Gr93}. It was answered negatively in 1998 by McMullen \cite{McM98} and also independently by Burago and Kleiner \cite{BK98}. Their results imply, in particular, that there are separated nets in $\R^d$ which are not a bounded displacement of $\Z^d$, not even after rescaling. 

In this paper we deal with separated nets which are obtained from tilings of Euclidean spaces. When a tiling of $\R^d$ is given, by placing one point in each tile, and keeping the minimal distance property, one gets a separated net. Since we are studying the equivalence classes under bounded displacement, the positions of the points in the tiles does not matter. In particular a tiling of $\R^d$ gives rise to a separated net; more precisely, an equivalence class of nets, in both of the above senses. Our main objective is to prove the two following theorems:

\begin{thm}\label{main1}
Any separated net that corresponds to a primitive substitution tiling of $\R^2$ is biLipschitz to $\Z^2$.
\end{thm}

\begin{definition}
Let $H$ be a primitive substitution in $\R^d$ and denote by $\lambda_2$ an eigenvalue of $A_H$ which is second in absolute value (see Definitions \ref{SubRule}, \ref{SubsTiling}, \ref{SubMatrix+Primitive}). If $\absolute{\lambda_2}<1$ we say that $H$ is a \emph{Pisot substitution}. 
\end{definition}


\begin{thm}\label{main2}
Let $H$ be a Pisot substitution in $\R^d$. Then for every substitution tiling of $H$ there exists a constant $\beta$ and bounded displacement between the corresponding separated net $Y$ and $\beta\cdot\Z^d$. 
\end{thm}

Here we also answer a question of Burago and Kleiner (\cite{BK02}, p.2). 

\begin{cor}\label{Penrose}
Any separated net that is created from a Penrose Tiling is a bounded displacement of $\beta\cdot\Z^2$, for some $\beta>0$ (and in particular biLipschitz to $\Z^2$).
\end{cor}

The proofs of Theorem \ref{main1} and of Theorem \ref{main2} rely on a result of Burago and Kleiner \cite{BK02} and a result of Laczkovich \cite{L92} respectively. Both of these results deal with the difference between the number of tiles in a large (bounded) set $U$ and the area of $U$. We use the Perron Frobenius Theorem and some dynamical properties of the matrix of the tiling in order to get good estimates for the number of tiles in large sets of certain kind. Then, by using properties of substitution tilings, we fill $U$ with such sets, that get smaller and smaller near the boundary of $U$, and get an estimate for the number of tiles in $U$. 


\medskip{\bf Acknowledgements:} This research was supported by the Israel Science Foundation. This work is a part of the author's Master thesis under the supervision of Barak Weiss whose endless support and guidance
are deeply appreciated. The author also wishes to thank Bruce Kleiner for suggesting the problem about the Penrose Tiling. After \cite{S08} appeared on the web it was brought to the author's attention that there is another paper, \cite{DSS95}, where a sketch of proof for Corollary \ref{Penrose} is given. 

\section{Basic Definitions of Tilings}
We use standard definitions of tilings. Similar definitions can be found at \cite{GS87}, \cite{Ra99}, \cite{Ro04}.

A set $S\subseteq\R^d$ is a \emph{tile} if it is homeomorphic to a closed $d$-dimensional ball. A \emph{tiling} of a set $U\subseteq\R^d$ is a countable collection of tiles, with pairwise disjoint interiors, such that their union is equal to $U$. We say that two tiles are \emph{translation equivalent} if one is a translation of the other. Representatives of the equivalence classes are called \emph{prototiles}. A \emph{tiling space}, $X_\mathcal{T}$, is the set of all tilings of $\R^d$ by prototiles from $\mathcal{T}$. A tiling $P$ of a bounded set $U\subset\R^d$ is called a \emph{patch}. We call the set $U$ the \emph{support of $P$} and we denote it by $supp(P)$. We extend the equivalence relation from the last Definition to patches and denote by $\mathcal{T}^*$ the equivalence class representatives.

\subsection*{Substitution Tilings}
Let $\xi>1$ and let $\mathcal{T}=\{S_1,\ldots,S_k\}$ be a set of $d$-dimensional prototiles. 
\begin{definition}\label{SubRule}
A \emph{substitution} is a mapping $H:\mathcal{T}\to\xi^{-1}\mathcal{T}^*$ such that for every $i$ we have $supp(S_i)=supp(H(S_i))$. In other words it is a set of dissection rules that shows us how to divide the prototiles to other prototiles from $\mathcal{T}$ with a smaller scale. We extend $H$ to the set of all tiles (in a given tiling), to $\mathcal{T}^*$ and to any tiling $\tau\in X_\mathcal{T}$ by applying $H$ separately on every tile. The constant $\xi$ is called the \emph{inflation constant} of $H$.
\end{definition}

\begin{definition}\label{SubsTiling}
Let $H$ be a substitution defined on $\mathcal{T}$. Consider the following set of patches: \[\mathcal{P}=\left\{(\xi H)^m(T_i):m\in\N, i=1,\ldots,k\right\}.\] The \emph{substitution tiling space} $X_H$ is the set of all tilings of $\R^d$ that for every patch $P$ in them there is a patch $P'\in\mathcal{P}$ such that $P$ is a sub-patch of $P'$. Every tiling $\tau\in X_H$ is a \emph{substitution tiling} of $H$. 
\end{definition}

\begin{prop}
If $H$ is a primitive substitution then $X_H\neq\emptyset$ and for every $\tau\in X_H$ and for every $m\in\N$ there exists a tiling $\tau_m\in X_H$ that satisfies $(\xi H)^m\tau_m=\tau$.   
\end{prop} 

\begin{proof}
See \cite{Ro04}.
\end{proof}

The construction of substitution tilings is explained with more details in \cite{Ro04}. We denote by $H^{(-m)}(\tau)$ a tiling $\tau '$ that satisfies $(\xi H)^m\tau '=\tau$. 

\subsection*{Matrices of Substitution}
\begin{Definition}
A matrix $A$ is called \emph{positive}, and denoted $A>0$, if all its entries are positive. $A$ is called \emph{nonnegative}, and denoted $A\ge 0$, if the entries of $A$ are nonnegative. A is called \emph{primitive} if there exists an $m\in \N$ such that $A^m>0$.
\end{Definition}

\begin{definition}
For a substitution $H$, the \emph{representative matrix} of $H$ is a $k\times k$ matrix $B_H=(b_{ij})$, where $b_{ij}$ is the number of prototiles which are translation equivalent to $S_i$ in $\xi H(S_j)$. We say that $H$ is \emph{primitive} if $B_H$ is primitive.
\end{definition}   

The matrix $B_H$ can be very large sometimes and does not describe exactly what we need here. Consider the following equivalence relation on prototiles: $S_i\sim S_j$ if there exists an isometry $O$ such that $S_i=O(S_j)$ and $H(S_i)=O(H(S_j))$ (this is actually a condition on the representatives, and it is obvious that it is well defined). We call the representatives of the equivalence classes \emph{basic tiles}. By this definition, we can also think of $H$ as a dissection rule on the basic tiles and extend it to tiles, patches and tilings as before. 

\begin{definition}\label{SubMatrix+Primitive}
Denote by $\{T_1,\ldots,T_n\}, (n\le k)$ the set of the basic tiles.  
Define the \emph{substitution matrix} of $H$ to be an $n\times n$ matrix, $A_H=(a_{ij})$, where $a_{ij}$ is the number of basic tiles in $\xi H(T_j)$ which are equivalent to $T_i$. 
\end{definition}   

\begin{example}
Let $H$ be the substitution of the Penrose Tiling. There are 20 different prototiles (rotations and reflections are not allowed):
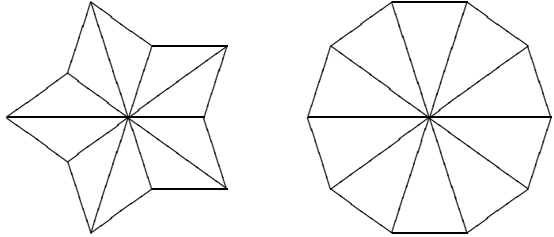
\[\xygraph{
!{<0cm,0cm>;<1cm,0cm>:<0cm,1cm>::}
!{(4,0) }*{}="a"
!{(5.618,0) }*{}="a1"
!{(5.309,0.951) }*{}="a2"
!{(4.5,1.534) }*{}="a3"
!{(3.5,1.534) }*{}="a4"
!{(2.691,0.951) }*{}="a5"
!{(2.382,0) }*{}="a6"
!{(2.691,-0.951) }*{}="a7"
!{(3.5,-1.534) }*{}="a8"
!{(4.5,-1.534) }*{}="a9"
!{(5.309,-0.951) }*{}="a10"
!{(0,0) }*{}="b"
!{(1,0) }*{}="b1"
!{(1.309,0.951) }*{}="b2"
!{(0.309,0.951) }*{}="b3"
!{(-0.5,1.534) }*{}="b4"
!{(-0.809,0.587) }*{}="b5"
!{(-1.618,0) }*{}="b6"
!{(-0.809,-0.587) }*{}="b7"
!{(-0.5,-1.534) }*{}="b8"
!{(0.309,-0.951) }*{}="b9"
!{(1.309,-0.951) }*{}="b10"
"a"-"a1" "a"-"a2" "a"-"a3" "a"-"a4" "a"-"a5"
"a"-"a6" "a"-"a7" "a"-"a8" "a"-"a9" "a"-"a10"
"a10"-"a1" "a1"-"a2" "a2"-"a3" "a3"-"a4" "a4"-"a5"
"a5"-"a6" "a6"-"a7" "a7"-"a8" "a8"-"a9" "a9"-"a10"
"b"-"b1" "b"-"b2" "b"-"b3" "b"-"b4" "b"-"b5"
"b"-"b6" "b"-"b7" "b"-"b8" "b"-"b9" "b"-"b10"
"b10"-"b1" "b1"-"b2" "b2"-"b3" "b3"-"b4" "b4"-"b5"
"b5"-"b6" "b6"-"b7" "b7"-"b8" "b8"-"b9" "b9"-"b10"
}\] 
Then $B_H$ is a $20\times 20$ matrix. On the other hand, there are only two different basic tiles, with the following dissection rule:
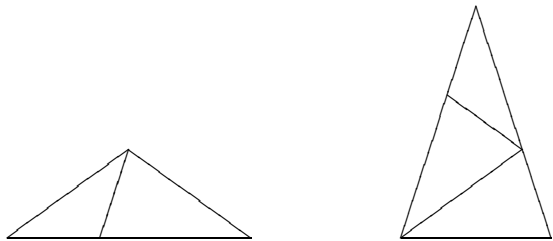
\[\xygraph{
!{<0cm,0cm>;<1cm,0cm>:<0cm,1cm>::}
!{(1,0) }*{}="a1"
!{(3,0) }*{}="a2"
!{(2,3.078) }*{}="a3"
!{(2.618,1.175) }*{}="sa1"
!{(1.618,1.902) }*{}="sa2"
!{(-1,0) }*{}="b1"
!{(-4.236,0) }*{}="b2"
!{(-2.618,1.175) }*{}="b3"
!{(-3,0) }*{}="sb"
"a1"-"a2" "a2"-"a3" "a3"-"a1" "a1"-"sa1" "sa1"-"sa2"
"b1"-"b2" "b2"-"b3" "b3"-"b1" "b3"-"sb"
}\] 
then $A_H=\begin{bmatrix} 2 & 1 \\ 1 & 1 \end{bmatrix}$, a $2\times 2$ matrix.  
\end{example}

Denote by $e_i$ the $i$'th element of the standard basis of $\R^n$ (or $\R^k$). Then $A_H(e_i)$ is the $i$'th column of $A_H$ (and the same for $B_H$ in $\R^k$). Thus, if $e_i$ represents one tile of type $i$, multiplying the vector $e_i$ by these matrices gives us a vector that represents the number of basic tiles (prototiles) of each kind obtained after applying $H$ on the corresponding tile. By linearity, this is true for any vector in $\R^n$ (or $\R^k$). Denote by $\pi:\R^k\to\R^n$ the quotient map that defines the relation $\sim$. Then by the definition of $\sim$, $\ker(\pi)$ is $B_H$-invariant and the following diagram commutes
\begin{diagram*}{
\R^k \ar[d]^{\pi} \ar[r]^{B_H} & \R^k \ar[d]^{\pi} \\
\R^n              \ar[r]^{A_H} & \R^n              &
}
\end{diagram*} 
Hence, it is easy to verify that the eigenvalues of $B_H$, that has an eigenvector $v\notin\ker(\pi)$, are also eigenvalues of $A_H$. 

\section{Properties of Substitution}
We denote by $\R^n_+$ the set of all nonnegative vectors in $\R^n$. For a finite set $P$ we denote by $\#P$ the number of elements of $P$. We also use the notations $\mu_d(\cdot), \norm{\cdot}_\infty$ and $\norm{\cdot}_2$ for $d$-dimensional Lebesgue measure, the max norm and the Euclidean norm in $\R^n$ respectively.

\begin{Definition}
Let $\lambda_1,\ldots,\lambda_n$ be the eigenvalues of a matrix $A$. The \emph{spectral radius} of $A$ is  $\rho(A)=\max_i\{\absolute{\lambda_i}\}$. 
For a nonnegative primitive matrix $A$, an eigenvalue $\lambda$ that satisfies $\absolute{\lambda}=\rho(A)$ is called \emph{a Perron Frobenius eigenvalue}.
\end{Definition}

Let $H$ be a primitive substitution with an inflation constant $\xi$ and let $\tau_0\in X_H$. Denote by $\{T_1,\ldots,T_n\}$ the set of $d$-dimensional basic tiles. We denote by $\lambda_1$ the Perron Frobenius eigenvalue of $A_H$ and let $\lambda_2$ be an eigenvalue which is second in absolute value. For every $m\in\N$ we denote $\tau_m=H^{(-m)}(\tau_0)$, these are substitution tilings with basic tiles $\{\xi^mT_1,\ldots,\xi^mT_n\}$. Then for every patch $P$ of $\tau_m$, $H^m(P)$ is a patch of $\tau_0$ and $supp(H^m(P))=supp(P)$. We denote by $t_i^{(m)}$ the number of tiles from $H^m(P)$ which are equivalent to $T_i$. 

Our main objective in this section is to prove the following proposition:
\begin{prop}\label{area_eval}
Let $H$ be a primitive substitution, then there are constants $a_1, a_2$ and $C_2$, that depend only on $H$, such that for every $0<\epsilon<\lambda_1-\absolute{\lambda_2}$ and for any $\tau_0\in X_H$ there exists an $N$ such that for every $m\ge N$ and a patch $P\in\tau_m$ we have
\begin{equation}\label{vol_eval_T}\begin{split}
t^{(m)}_1(a_1-C_2\delta^m)\le\#P\le t^{(m)}_1(a_1+C_2\delta^m) \\
t^{(m)}_1(a_2-C_2\delta^m)\le \mu_d(V) \le t^{(m)}_1(a_2+C_2\delta^m),
\end{split}\end{equation}
where $V=supp(P)$ and 
\begin{equation}\label{delta}
\delta=\frac{\absolute{\lambda_2}+\epsilon}{\lambda_1}<1.
\end{equation} 
\end{prop}

We start with some notions of matrix theory. See also \cite{BP79}, \cite{H07}.

\begin{prop}
Let $H$ be a primitive substitution with an inflation constant $\xi$ and let $\{T_1,\ldots,T_n\}$ be the set of $d$-dimensional basic tiles, then
\begin{description}\label{PF}
    \item[(a)] $\rho(A_H)$ is an eigenvalue of $A_H$, with an algebraic multiplicity one, and the associate eigenvector is positive.
    \item[(b)] \label{B} If $v>0$ is an eigenvector of $A_H$ then $v$ corresponds to $\rho(A_H)$.
     \item[(c)] Denote by $(v_1,\ldots,v_n)$ a basis of generalized eigenvectors of $A_H$ with $\norm{v_i}_2=1$ for all $i$, where $v_1$ is the eigenvector that corresponds to $\rho(A_H)$. Then $sp\{v_2,\ldots,v_n\}\cap\R^n_+=\{0\}$.
     \item[(d)] In addition, $\xi^d=\rho(A_H)$, and in particular $\rho(A_H)>1$.   
\end{description}
\end{prop}

\begin{proof}
\begin{description}
    \item[(a)] This is the well known Perron Frobenius Theorem, also see \cite{Q87}, p.91. 
    \item[(b)] See \cite{BP79} Theorem 2.1.4.
    \item[(c)] This is a well known fact, see \cite{H07}, 26.1.4.(d).    
    \item[(d)] Consider the vector $v=\begin{bmatrix} s_1 \\ \vdots \\ s_n  \end{bmatrix}$ where $s_i$ is the area of $T_i$. Since the substitution divides every $T_i$ to smaller basic tiles, by using Definition \ref{SubMatrix+Primitive} one can easily show that $\xi^d v=A_H^tv$, where $A_H^t$ is the transpose matrix of $A_H$. By {\bf (b)} the proof is complete. 
    
\end{description}
\end{proof}

Let $A\ge 0$ be a primitive matrix. We use the notations $\lambda_1$ and $\lambda_2$ as before. We also denote by $v_1$ the eigenvector of $\lambda_1$ with $\norm{v_1}_2=1$ and by $W$ the generalized eigenspace of the other eigenvalues. For a vector $u\in\R^n$ we write $u=\beta_1(u)v_1+\beta_2(u)w_u$ where $w_u\in W$ with $\norm{w_u}_2=1$.      

\begin{prop}\label{convW1}
For a primitive matrix $A\ge 0$ there exists a constant $C>0$, that depends only on $A$, such that for every $\epsilon>0$ there exists an $N$ such that for every $m\ge N$ and a vector $u\in\R^n_+$ we have 
\begin{equation}\label{vector_of_tiles}
\norm{\frac{A^mu}{\beta_1(u)\lambda_1^m}-v_1}_\infty \le 
\delta^m,
\end{equation}
where $\delta$ as in (\ref{delta}).
\end{prop}

\begin{proof}
We denote $K=\{v\in\R^n_+:\norm{v}_2=1\}$ and consider the restrictions of $\beta_1$ and $\beta_2$ to $K$. 
By Proposition \ref{PF} (c) we have $\beta_1(u)>0$ for every $u\in K$. By the compactness of $K$ we denote 
\[\alpha_1=\min_{u\in K}\{\beta_1(u)\}>0 \quad\mbox{and}\quad \alpha_2=\max_{u\in K}\{\absolute{\beta_2(u)}\}>0.\]
Denote $C'=\frac{\alpha_2}{\alpha_1}$, then $C'>0$ and for every $u\in\R^n_+$ we have: 
\[\frac{\absolute{\beta_2(u)}}{\beta_1(u)}= \frac{\norm{u}_2\cdot\absolute{\beta_2(\frac{u}{\norm{u}_2})}}{\norm{u}_2\cdot \beta_1(\frac{u}{\norm{u}_2})} \le\frac{\alpha_2}{\alpha_1}=C'.\]
Notice that there is a constant $a>0$ such that for every $m\in\N$ and $w\in W$ with $\norm{w}_2=1$ we have $\norm{A^m(w)}_\infty\le a(\absolute{\lambda_2}+\epsilon)^m$. Then 
\[\norm{\frac{A^mu}{\beta_1(u)\lambda_1^m}-v_1}_\infty= \norm{\frac{A^m(\beta_1(u)v_1+\beta_2(u)w_u)-\beta_1(u)\lambda_1^mv_1}{\beta_1(u)\lambda_1^m}}_\infty\]
\[=\norm{\frac{A^m(\beta_2(u)w_u)}{\beta_1(u)\lambda_1^m}}_\infty= \frac{\absolute{\beta_2(u)}\cdot\norm{A^m(w_u)}_\infty}{\beta_1(u)\lambda_1^m}\le C'\cdot\frac{a(\absolute{\lambda_2}+\epsilon)^m}{\lambda_1^m},\]
which completes the proof for $C=C'a$.  
\end{proof}

\begin{prop}\label{tilesratio}
There are constants $C_1,c_2,\ldots,c_n>0$, that depend only on $H$, such that for every $0<\epsilon<\lambda_1-\absolute{\lambda_2}$ there exists an $N_1$ such that every $m\ge N_1$ satisfies
\begin{equation}\label{conv_to_ratio}
\absolute{\frac{t_i^{(m)}}{t_1^{(m)}}-c_i}\le C_1\delta^m
\end{equation}
for every patch $P$ of $\tau_m$, where $\delta$ as in (\ref{delta}).  
\end{prop}

\begin{proof}
By Proposition \ref{PF}, $\lambda_1>1$ and it has an associated eigenvector $v_1=\begin{bmatrix}1 \\ c_2 \\ \vdots \\ c_n \end{bmatrix}$ ($c_1=1$) with $c_i>0$ for $i=2,\ldots,n$. Denote $v_1'=\frac{v_1}{\norm{v_1}_2}$. Fix $0<\epsilon<\lambda_1-\absolute{\lambda_2}$, then there is an $N$ such that for every $m\ge N$ and $u\in\R^n_+$ we have (\ref{vector_of_tiles}) with $v_1'$ instead of $v_1$. 

We pick $N_1\ge N$ such that every $m\ge N_1$ satisfies
\begin{equation}\label{half}
C\norm{v_1}_2\delta^m\le \frac{1}{2}.
\end{equation}

For an arbitrary $m\ge N_1$ and a patch $P$ of $\tau_m$, consider the vector $u=\begin{bmatrix} t_1^{(0)} \\ \vdots \\ t_n^{(0)} \end{bmatrix}$, where $t_i^{(0)}$ is the number of tiles from $P$ which are equivalent to $\xi^mT_i$.    Obviously $u\in\R^n_+\minus\{0\}$, then $u$ satisfies (\ref{vector_of_tiles}) with $v_1'$ instead of $v_1$.
Hence
\[\norm{\frac{\norm{v_1}_2\cdot A_H^mu}{\alpha_1(u)\lambda_1^m}-v_1}_\infty \le C\norm{v_1}_2\delta^m.\]
If we denote $\frac{\norm{v_1}_2\cdot A_H^mu}{\alpha_1(u)\lambda_1^m}= 
\begin{bmatrix}b_1^{(m)} \\ \vdots \\ b_n^{(m)} \end{bmatrix}$ then for $i=2,\ldots,n$ we have 
\begin{equation}\label{close_in_coordinate}
\absolute{b_i^{(m)}-c_i}\le C\norm{v_1}_2\delta^m \quad\mbox{and}\quad \absolute{b_1^{(m)}-1}\le C\norm{v_1}_2\delta^m.
\end{equation}
In particular, by (\ref{half}), $\frac{1}{2}\le b_1^{(m)}\le 1\frac{1}{2}$, for every $m\ge N_1$. 
Notice that by the definitions of $u$, $P$ and $A_H$ we have  
\begin{equation}\label{A_H^mu}
A_H^mu=\begin{bmatrix} t^{(m)}_1 \\ \vdots \\ t^{(m)}_n \end{bmatrix}.
\end{equation}  
Therefore for $i=2,\ldots,n$ we have \[\absolute{\frac{t^{(m)}_i}{t^{(m)}_1}-c_i}\stackrel{(\ref{A_H^mu})}= \absolute{\frac{b_i^{(m)}}{b_1^{(m)}}-c_i}\le \frac{1}{\absolute{b_1^{(m)}}}\left(\absolute{b_i^{(m)}-c_i}+\absolute{c_i-b_1^{(m)}c_i}\right)\]
\[\stackrel{(\ref{half}), (\ref{close_in_coordinate})}\le 2C\norm{v_1}_2(1+c_i)\cdot\delta^m \le C_1\cdot\delta^m,\]
where $C_1=2C\norm{v_1}_2(1+\max_i\{c_i\})$, as required. 
\end{proof}

We denote by $s_1,\ldots,s_n$ the areas of $\{T_1,\ldots,T_n\}$ respectively. Define
\begin{equation}\label{alpha}
a_1=\sum_{i=1}^nc_i\quad,\quad a_2=\sum_{i=1}^nc_is_i \quad\mbox{ and }\quad
\alpha=\frac{a_1}{a_2}
\end{equation}

\begin{proof}[Proof of Proposition \ref{area_eval}]
Let $0<\epsilon<\lambda_1-\absolute{\lambda_2}$. By Proposition \ref{tilesratio} there exists an $N=N_1$ such that (\ref{conv_to_ratio}) holds for every $m\ge N$ and a patch $P$ in $\tau_m$, for some constant $C_1$. 
Then for $i=2,\ldots,n$ we have 
\[t^{(m)}_1(c_i-C_1\cdot\delta^m) \le t^{(m)}_i\le t^{(m)}_1(c_i+C_1\cdot\delta^m)\] and
\begin{align*}
t^{(m)}_1\left(c_i s_i-C_1\cdot\delta^m s_i\right) \le&t^{(m)}_i s_i\le t^{(m)}_1\left(c_is_i+C_1\cdot\delta^m s_i\right).
\end{align*}
Therefore
\begin{align*}
t^{(m)}_1\left(\sum_{i=1}^nc_i-nC_1\cdot\delta^m\right) 
\le&\sum_{i=1}^n t^{(m)}_i \le t^{(m)}_1\left(\sum_{i=1}^nc_i+nC_1\cdot\delta^{m}\right)
\\
t^{(m)}_1\left(\sum_{i=1}^n c_i s_i-C_1\cdot\delta^m \sum_{i=1}^n s_i\right) \le&\sum_{i=1}^n t^{(m)}_i s_i\le t^{(m)}_1\left(\sum_{i=1}^nc_is_i+C_1\cdot\delta^m\sum_{i=1}^ns_i\right).
\end{align*}
Thus, according to (\ref{alpha}), for $C_2=\max\{C_1n,C_1\sum_{i=1}^ns_i\}$ we get (\ref{vol_eval_T}) as required.
\end{proof}

\section{The Main Results}

We prove Theorem \ref{main1} by showing that substitution tilings create separated nets that satisfy the conditions of the following theorem: 
\begin{thm}[Burago and Kleiner \cite{BK02}]
Let Y be a separated net in $\R^2$. For a real number $\alpha>0$ and a square $B$ with integer coordinates define:
\begin{align*} \label{E} e_{\alpha}(B)=&\max\bigg\{\frac{\alpha\cdot\mu_2(B)}{\#(B\cap Y)}, \frac{\#(B\cap Y)}{\alpha\cdot\mu_2(B)} \bigg\} \\ 
E_{\alpha}(2^i)=&\sup\big\{e_{\alpha}(B) : B\quad \mbox{as above with an edge of length } 2^i \big\}.
\end{align*}  
If there exists an $\alpha>0$ such that the product $\prod_{j=1}^{\infty} E_{\alpha}(2^j)$ converges, then $Y$ is biLipschitz to $\mathbb{Z}^2$.
\end{thm}

\begin{proof}[Proof of Theorem \ref{main1}]
It suffices to show that there are constants $C_1,k_1>0$ and $\omega<1$ such that for every square $B$ with an edge of length $2^j=k\ge k_1$ we have 

\begin{equation}\label{BK}
\frac{\absolute{\alpha\cdot\mu_2(B)-\#(B\cap Y)}}{\#(B\cap Y)} \le C_1\cdot \omega^j\quad\mbox{ and }\quad \frac{\absolute{\alpha\cdot\mu_2(B)-\#(B\cap Y)}}{\alpha\cdot\mu_2(B)} \le C_1\cdot \omega^j.
\end{equation} 
Then we get, for all large enough $j$, 
\[\absolute{E_\alpha(2^j)-1}\le C_1\cdot \omega^j,\]
which implies the convergence of the product.

By Proposition \ref{PF} we can pick an $\epsilon>0$ such that $\lambda_1>\absolute{\lambda_2}+\epsilon$.  By Proposition \ref{area_eval} there is an $N_1$ such that for every $m\ge N_1$ and a patch $P$ in $\tau_m$, (\ref{vol_eval_T}) holds. Let $N\ge N_1$ such that for every $m\ge N$ we have 
\begin{equation}\label{new_N_1.3}
C_2\cdot\delta^m\le\frac{1}{2}\min\{a_1,a_2\}\quad (a_1, a_2 \mbox{ as in } (\ref{alpha}), \delta\mbox{ as in } (\ref{delta})).
\end{equation} 
We pick $k'=\xi^{2N}$. 

Let $B$ be an arbitrary square in $\R^2$ with an edge of length $k\ge k'$. Let $m\in\N$ such that $\xi^{2m}\le k<\xi^{2m+2}$, then $m\ge N$. Consider the patch 
\[P=\{T\in\tau_m : T\subseteq B\},\]
then $P$ satisfies (\ref{vol_eval_T}), where $V=supp(P)$ as before.

Let $R$ and $r$ be constants such that every ball of diameter $R$ contains a tile of $\tau_0$ and every tile of $\tau_0$ contains a cube of area $r$. 
From the definition of $P$, for every $x\in B$ that satisfies $d(x,\partial B)\ge R\cdot\xi^{m}$, the tile of $\tau_m$ that covers $x$ must be in $P$. Then $V$ contains a square with an edge of length $k-2R\cdot\xi^m$. Since $R\cdot\xi^m\le R\cdot\sqrt{k}$, $V$ contains a square with an edge of length $k-2R\cdot\sqrt{k}$. If so, there is a $k''$ such that for every $k\ge k''$ we have $\mu_2(V)\ge\frac{1}{2}k^2$. Then, by (\ref{vol_eval_T}), there is a constant $b_1>0$ such that for every $k\ge k''$ we have 
\begin{equation}\label{t_1^m}
b_1\cdot k^2\le t^{(m)}_1.
\end{equation} 

Define $k_1=\max\{k',k''\}$. Consider squares $B$ with an edge of length $k\ge k_1$. We want to estimate $\#(B\cap Y)$ and $\mu_2(B)$. Define the following patch of $\tau_0$:
\[P_1=\{T\in\tau_0:T\cap B\neq\emptyset\} \qquad V_1=supp(P_1).\] 

A similar explanation to the one above gives the estimate
\[V_1\minus V\subseteq\{x: d(x,\partial B)\le R\cdot\xi^{m}\}\subseteq\{x: d(x,\partial B)\le R\cdot\sqrt{k}\}.\] 
Then $\mu_2(V_1\minus V)\le 4R\cdot k\sqrt{k}$ and so $\#((V_1\minus V)\cap Y)\le\frac{4R\cdot k\sqrt{k}}{r}$. Therefore 

\begin{align*}
\#P\le \#(B\cap Y)&\le \#P+\#((V_1\minus V)\cap Y), \\ \\
\mu_2(V)\le \mu_2(B)&\le\mu_2(V)+\mu_2(V_1\minus V).
\end{align*}
Hence, by (\ref{vol_eval_T}) 
\begin{align*}
t^{(m)}_1(a_1-C_2\delta^m)\le \#(B\cap Y)\le t^{(m)}_1(a_1+C_2c\delta^m)+\frac{4R\cdot k\sqrt{k}}{r}, \\ 
t^{(m)}_1(a_2-C_2\delta^m)\le \mu_2(B) \le t^{(m)}_1(a_2+C_2\delta^m)+4R\cdot k\sqrt{k}.
\end{align*}
Therefore, for $\alpha$ as in (\ref{alpha}) we have
\begin{align*}
\frac{\alpha\cdot\mu_2(B)-\#(B\cap Y)}{\#(B\cap Y)}
\le\frac{\alpha(t^{(m)}_1(a_2+C_2\delta^m)+4R\cdot k\sqrt{k} )-t^{(m)}_1(a_1-C_2\delta^m)}{t^{(m)}_1(a_1-C_2\delta^m)} \\ 
=\frac{t^{(m)}_1 C_2\delta^m(\alpha+1)+4\alpha R\cdot k\sqrt{k}}{t^{(m)}_1(a_1-C_2\delta^m)}
\stackrel{(\ref{t_1^m})}\le\delta^m\cdot\frac{C_2(\alpha+1)}{a_1-C_2\delta^m}+\frac{1}{\sqrt{k}}\cdot \frac{4\alpha R}{b_1(a_1-C_2\delta^m)}.
\end{align*}
In the same way we get similar inequalities  for $\frac{\alpha\cdot\mu_2(B)-\#(B\cap Y)}{\#(B\cap Y)}$ and then for $\frac{\absolute{\alpha\cdot\mu_2(B)-\#(B\cap Y)}}{\alpha\cdot\mu_2(B)}$. Hence, considering (\ref{new_N_1.3}), there is a constant $C_1'$ such that
\[\max\left\{\frac{\absolute{\alpha\cdot\mu_2(B)-\#(B\cap Y)}}{\alpha\cdot\mu_2(B)},\frac{\absolute{\alpha\cdot\mu_2(B)-\#(B\cap Y)}}{\#(B\cap Y)} \right\} \le\left(\delta^m+\frac{1}{\sqrt{k}}\right)\cdot C_1'.\]
Notice that $m$ was chosen in a way that $\xi^{2m}\le k$, thus $\xi^m\le \sqrt{k}$. Since we are looking on squares with $k=2^j$, we get that $m\le j\cdot\log_\xi \sqrt{2}$. Therefore, $\omega=\max\{\delta^{\log_\xi \sqrt{2}},1/\sqrt{2}\}$ satisfies the condition in (\ref{BK}), with $C_1=2C_1'$.
\end{proof}

We now turn to the proof of Theorem \ref{main2}. The proof relies on the following theorem:

\begin{thm}[Laczkovich \cite{L92}]
For a separated net $Y\subseteq\R^d$ and $\alpha>0$ the following statements are equivalent:\newline
(i) There is a positive constant $C$ such that for every finite union of unit cubes $U$ we have  
\begin{equation}\label{Laczkovich}
\absolute{\#({Y\cap U})-\alpha\mu_d(U)}\le C\cdot\mu_{d-1}(\partial U). \end{equation}  
(ii) There is a bounded displacement $\phi:Y\to\alpha^{-1/d}\Z^d$.
\end{thm}

For the proof of Theorem \ref{main2} we will need the two following lemmas:

\begin{lem}\label{one_tile}
There is a constant $C_3$ such that for every $0<\epsilon<\lambda_1-\absolute{\lambda_2}$ there exists an $N$ such that for every $m\ge N$ and a tile $T$ in $\tau_m$ we have
\begin{equation}\label{Laczkovich_T}
\absolute{\#(T\cap Y)-\alpha\mu_d(T)}\le C_3\cdot(\lambda_2+\epsilon)^m,
\qquad\mbox{where $\alpha$ as in (\ref{alpha}).}
\end{equation}

\end{lem}
\begin{proof}
We think of $T$ as a patch in $\tau_m$ and denote $P_0=H^m(T)$, the patch in $\tau_0$ with $supp(P_0)=T$. Then by Proposition \ref{area_eval} we have (\ref{vol_eval_T}) with $T$ instead of $P$, for every $m\ge N_1$. On the other hand, $T$ is equivalent to $\xi^mT_i$ for some $i\in\{1,\ldots,n\}$, then $\mu_d(T)=(\xi^d)^m\cdot s_i$. Then
\[t^{(m)}_1(a_2-C_2\delta^m)\le(\xi^d)^m\cdot s_i,\] 
which implies, for every $m$ which is greater than some $N_2$, 
\begin{equation}\label{t_1^(m)}
t^{(m)}_1\le\frac{(\xi^d)^m\cdot s_i}{a_2-C_2\delta^m}\le C_3'(\xi^d)^m \stackrel{(\ref{PF})(c)}=C_3'\lambda_1^m.
\end{equation}

According to (\ref{vol_eval_T}) we have
\begin{align*}
\#(T\cap Y)-\alpha\mu_d(T)\stackrel{(\ref{vol_eval_T})}\le t^{(m)}_1(a_1+C_2\delta^m)-\alpha t^{(m)}_1(a_2-C_2\delta^m) \stackrel{(\ref{alpha})}=t^{(m)}_1C_2\delta^m(1+\alpha),
\end{align*}
and in a similar way we get it for $\alpha\mu_d(T)-\#(T\cap Y)$. Then 
\[\absolute{\#(T\cap Y)-\alpha\mu_d(T)}\le t^{(m)}_1C_2\delta^m(1+\alpha) \stackrel{(\ref{t_1^(m)})}\le C_3'\lambda_1^m \left(\frac{\absolute{\lambda_2}+\epsilon}{\lambda_1}\right)^m C_2(1+\alpha).\]
All this is true for every $m\ge N$, where $N=\max\{N_1,N_2\}$. Then for $C_3=C_3'\cdot C_2(1+\alpha)$ we get the required inequality.
\end{proof}

\begin{lem}
There is a constant $C$, that depends only on the dimension $d$, such that for any $s>1$  
\begin{equation}\label{Lacz.inequality}
\mu_d\left(\{x\in U: d(x,\partial U)\le s \}\right)\le C\cdot s^d\cdot\mu_{d-1}(\partial U)
\end{equation} 
holds for any finite union of $d$-dimensional cubes $U$.
\end{lem}

\begin{proof}
This a direct result of Lemma 2.1 and Lemma 2.2 of \cite{L92}.
\end{proof}

\begin{proof}[Proof of Theorem \ref{main2}]
First we claim that it is sufficient to show inequality (\ref{Laczkovich}) for every set $U$ which is a finite union of cubes with an edge of length $k$, for some constant $k\in\N$. Then indeed, by rescaling the whole picture by a factor of $\frac{1}{k}$, we get (\ref{Laczkovich}) for the net $\frac{1}{k}\cdot Y$ with $\frac{1}{k}\cdot U$ instead of $U$, $\frac{\alpha}{k^d}$ instead of $\alpha$ and a different constant $C_1$. Since $\frac{1}{k}\cdot U$ is a finite union of unit cubes, by \cite{L92}, we get a 
bounded displacement ${\Phi}':\frac{1}{k}\cdot Y\to\alpha^{1/d}\frac{1}{k}\cdot\Z^d$, which implies the existence of the required $\Phi$.

We pick the constant $\alpha$ as in (\ref{alpha}).
Since $H$ is a Pisot substitution, we fix $\epsilon>0$ such that $\absolute{\lambda_2}+\epsilon<1$ (By Proposition \ref{PF} $\lambda_1=\xi^d >1$). Then let $N$ be such that (\ref{Laczkovich_T}) holds for every tile $T\in\tau_m$ where $m\ge N$. Let $R$ and $r$ be constants such that every ball of diameter $R$ contains a tile of $\tau_0$ and every tile of $\tau_0$ contains a cube of area $r$. 
Let $U$ be a finite union of cubes in $\R^d$ with an edge of length $k=\ceil{R\cdot\xi^N}$. Let $m$ be the maximal integer such that $U$ contains a tile of $\tau_m$. Then by the definition of $k$ we have $m\ge N$. Define the following sequence of patches: \[P_m=\left\{T\in\tau_m:T\subseteq U\right\}\]
and for decreasing $l=m-1,\ldots,N$
\[P_l=\left\{T\in\tau_l:int(T)\subseteq U\minus\bigcup_{j=l+1}^m V_j\right\},\]
where $V_l=supp(P_l)$ and $int(T)$ is the interior of $T$. 

Define $V_{\partial}=U\minus\bigcup_{j=N}^m V_j$, then we get a partition of $U$ to layers that intersect only at their boundaries: \begin{equation}\label{U_union_of_cubes}
U=\left(\bigcup_{l=N}^m V_l\right)\cup V_{\partial},
\end{equation} which implies
\begin{equation*}
\mu_d(U)=\left(\sum_{l=N}^m\mu_d(V_l)\right)+\mu_d(V_{\partial}).
\end{equation*}  
 
We now estimate $\#P_l$.
Notice that for every $x\in U$, if $d(x,\partial U)\ge R\cdot\xi^l$ then any ball of diameter $R\cdot\xi^l$ that contains $x$ is contained in $U$. Then the tile of $\tau_l$ that contains $x$ is contained in $U$. In particular we get it for $l=m+1$. But since no tile of $\tau_{m+1}$ is contained in $U$, we deduce that $d(x,\partial U)<R\cdot\xi^{m+1}$ for every $x\in U$. Therefore \[\mu_d(U)\le \mu_d\left(\left\{x\in U:d(x,\partial U)<R\cdot\xi^{m+1}\right\}\right) \stackrel{(\ref{Lacz.inequality})}\le C\cdot (R\xi^{m+1})^d\mu_{d-1}(\partial U).\]  
Since every tile of $\tau_m$ contain a cube of area $r(\xi^d)^m$ we have
\[\#P_m\le\frac{C(R\xi^{m+1})^d\mu_{d-1}(\partial U)}{r(\xi^d)^m}= \frac{C(R\xi)^d\mu_{d-1}(\partial U)}{r}.\]    
In a similar way, for every $l$, if $d(x,\partial U)\ge R\cdot\xi^{l+1}$ then the tile of $\tau_{l+1}$ that covers $x$ is in $P_{l+1}$, thus $x\notin V_l$. Hence for every $x\in V_l$ we have $d(x,\partial U)<R\cdot\xi^{l+1}$ and so \[\mu_d(V_l)\le \mu_d\left(\left\{x\in U:d(x,\partial U)<R\cdot\xi^{l+1}\right\}\right) \stackrel{(\ref{Lacz.inequality})}\le C\cdot (R\xi^{l+1})^d\mu_{d-1}(\partial U),\] 
which implies
\begin{equation}\label{P_l} 
\#P_l\le\frac{C(R\xi^{l+1})^d\mu_{d-1}(\partial U)}{r(\xi^d)^l}= \frac{C(R\xi)^d\mu_{d-1}(\partial U)}{r}.
\end{equation} 
Similarly we get 
\begin{equation}\begin{split}\label{V_partial}
\mu_d(V_\partial)\le C(R\xi^N)^d\mu_{d-1}(\partial U) \\
\#(V_\partial\cap Y)\le \frac{C(R\xi^N)^d\mu_{d-1}(\partial U)}{r}.
\end{split}\end{equation}
If we denote by $T^{(l)}$ tiles of $\tau_l$, then for every $l\ge N$ we have 
\begin{gather*}
\absolute{\#(V_l\cap Y)-\alpha\mu_d(V_l)}= \absolute{\sum_{T^{(l)}\subseteq V_l}\#(T^{(l)}\cap Y)-\alpha\cdot\sum_{T^{(l)}\subseteq V_l}\mu_d(T^{(l)})} \\
\le\sum_{T^{(l)}\in P_l}\absolute{\#(T^{(l)}\cap Y)-\alpha\mu_d(T^{(l)})} \stackrel{(\ref{Laczkovich_T}),(\ref{P_l})}\le\frac{C(R\xi)^d\mu_{d-1}(\partial U)}{r}\cdot C_3\cdot(\absolute{\lambda_2}+\epsilon)^l \\
\le C_4\cdot(\absolute{\lambda_2}+\epsilon)^l \cdot\mu_{d-1}(\partial U),
\end{gather*}
where $C_4=\frac{C\cdot C_3\cdot(R\xi)^d}{r}$. Therefore, according to (\ref{V_partial}), we denote $C_5=\max\left\{\alpha\cdot C(R\xi^N)^d,\frac{(R\xi^N)^d}{r}\right\}$ and get
\begin{gather*}
\absolute{\#(U\cap Y)-\alpha\mu_d(U)}\stackrel{(\ref{U_union_of_cubes})} \le\left[\sum_{l=N}^m\absolute{\#(V_l\cap Y)-\alpha\mu_d(V_l)}\right]+ \absolute{\#(V_\partial\cap Y)-\alpha\mu_d(V_\partial)} \\
\le\left[\sum_{l=N}^m C_4\cdot(\absolute{\lambda_2}+\epsilon)^l\cdot\mu_{d-1}(\partial U) \right]+C_5\cdot\mu_{d-1}(\partial U)\le C_1\cdot\mu_{d-1}(\partial U),
\end{gather*}
where $C_1=C_4\left(\sum_{l=1}^\infty(\absolute{\lambda_2}+\epsilon)^l\right)+C_5$.
\end{proof}

\begin{proof}[Proof of Corollary \ref{Penrose}]
If we denote by $H$ the substitution of the Penrose Tiling then $A_H=\begin{bmatrix} 2 & 1 \\ 1 & 1 \end{bmatrix}$. In this case we have 
$\lambda_2=\frac{3-\sqrt{5}}{2}<1$. By Theorem \ref{main2} the proof is complete. 
\end{proof}


\begin{thebibliography}{100000}
\bibitem[BP79]{BP79} A. Berman and R. J. Plemmons, \underline{Nonnegative matrices in the} \underline{mathematical sciences}, Academic Press, New York, 1979.
\bibitem[BK98]{BK98} D. Burago and B. Kleiner, \emph{Separated nets in Euclidean space and Jacobians of biLipschitz map}, Geom. Func. Anal. 8 (1998), no.2, 273-282.
\bibitem[BK02]{BK02} D. Burago and B. Kleiner, \emph{Rectifying separated nets}, Geom. Func. Anal. Vol.12 (2002) 80-92.
\bibitem[DSS95]{DSS95} W. A. Deuber, M. Simonovits and V. T. Sos, \emph{A note on paradoxical metric spaces}, Studia Sci.Hung.Math. 30 (1995), 17--23. 
\bibitem[Gr93]{Gr93} M. Gromov, \underline{Asymptotic invariants of infinite groups}, Geometric Group Theory Vol.II (G. Niblo and M. Roller eds.), London Math. Soc. Lecture Notes , 182, Cambridge Univ. Press  (1993).
\bibitem[GS87]{GS87} Branko Grubaum and G.C. Shephard, \underline{Tilings and patterns}, W.H.Freeman and Company, New York, 1987.
\bibitem[H07]{H07} L. Hogben (Ed.), \underline{Handbook of linear algebra}, Chapman and Hall/CRC Press, 2007. MR2279160 (2007j:15001). 
\bibitem[L92]{L92} M. Laczkovich, \emph{Uniformly spread discrete sets in $\R^d$}, J. London Math. Soc. (2) 46 (1992) 39-57.
\bibitem[McM98]{McM98} C. T. McMullen, \emph{Lipschitz maps and nets in Euclidean space}, Geom. Func. Anal. 8 (1998), no.2, 273-282.
\bibitem[Q87]{Q87} M. Queffelec, \emph{Substitution dynamical systems-Spectral analysis}, Lecture notes in mathematics, vol.1294, Springer-Verlag, Berlin, 1987. 
\bibitem[Ra99]{Ra99} C. Radin, \underline{Miles of tiles}, Amer. Math. Soc., Providence, RI (1999).
\bibitem[Ro04]{Ro04} E. A. Robinson, Jr. \emph{Symbolic dynamics and tilings of $\R^n$}, Proc. Sympos. Appl. Math. Vol.60 (2004), 81-119.
\bibitem[S08]{S08} Y. Solomon, \emph{The net created from the Penrose Tiling
is biLipschitz to the integer lattice}, arXiv:0711.3707v1 (2008).

\end{thebibliography}
\end{document}